\documentclass[a4paper, 12pt]{article}

\usepackage[unicode]{hyperref}
\usepackage{amsmath}
\usepackage{amssymb}
\usepackage[strict,style=british]{csquotes}
\usepackage{graphicx}
\usepackage{booktabs}
\usepackage{geometry}
\usepackage{authblk}
\usepackage[sort]{cite}
\usepackage{tablefootnote}
\usepackage{algorithm}
\usepackage{hyperref}
\usepackage{url}
\usepackage[noend]{algpseudocode}
\usepackage{multirow}
\usepackage{comment}
\usepackage[sort]{cite}
\usepackage{caption}
\usepackage{subcaption}
\usepackage{epstopdf}

\title{Direct simulation Monte Carlo for new regimes in aggregation-fragmentation kinetics}
\author[1]{A.~Kalinov}
\author[1]{A.~.I.~Osinskiy}
\author[2,3]{S.~A.~Matveev}
\author[4]{W.~Otieno}
\author[1,4]{N.~V.~Brilliantov}

\date{}
\affil[1]{Skolkovo Institute of Science and Technology, Moscow, Russia}
\affil[2]{Faculty of Computational Mathematics and Cybernetics, Lomonosov Moscow State University, Russia}
\affil[3]{Marchuk Institute of Numerical Mathematics of Russian Academy of Sciences, Moscow, Russia}
\affil[4]{University of Leicester, Leicester, UK}

\begin{document}
\maketitle{}

\begin{abstract}
    We revisit two basic Direct Simulation Monte Carlo Methods to model aggregation kinetics and extend them for aggregation processes with collisional fragmentation (shattering). We test the performance and accuracy of the extended methods and compare their performance with efficient deterministic finite-difference method applied to the same model. We validate the stochastic methods on the test problems and apply them to verify the existence of oscillating regimes in the aggregation-fragmentation kinetics recently detected in deterministic simulations. We confirm the emergence of steady oscillations of densities in such systems and prove the stability of  the oscillations with respect to fluctuations and noise. 
\end{abstract}

\textbf{Keywords ---} coagulation, aggregation-fragmentation kinetics, Smoluchowski equations, Direct Simulation Monte Carlo, steady oscillations.

\section{Introduction}

Aggregation kinetics -- the process where particles (elementary units) join together to form larger agglomerates --  has been successfully explored experimentally and analytically for more than a century and numerically for a couple of decades. The physical nature of an elementary unit may be very different, ranging from molecules and colloids to fibrils and dust grains \cite{Leyvraz,KrapivskyBook,Krapivsky_PRB1999,SpahnetalEPL:2004,BrilliantovSpahn2006,DominikTilens:1997,PNAS,BFP2018,prions}. Moreover, aggregation is observed in networks such as Internet or business/social networks \cite{CoagNetw}.   Self-assembly -- an important natural phenomenon -- is also associated with aggregation  \cite{SA08,Igor14,Erik17}.  The history of the field commenced in 1914 from the seminal work of  Smoluchowski \cite{smoluchowski1916drei}, where the aggregation of colloidal particles in solution was put in a mathematical framework. Namely, famous Smoluchowski equations for concentrations $n_s(t)$  of aggregates of size $s$ have  been formulated: 

\begin{equation}
\label{SM}
\frac{dn_s}{d t} = \frac{1}{2} \sum_{i=1}^{s-1} K_{i, s-i} n_i n_{s-i} - n_s \sum_{i=1}^{\infty} K_{i,s} n_i.
\end{equation}
Eqs. (\ref{SM}) is an infinite system of the first-order ODE. The first term on the r.h.s. of (\ref{SM}) gives the rate of the process when two particles of size $i$ and $(s-i)$ coalesce, giving rise to a particle of size $s$. The summation accounts for all such processes, which increase the concentration $n_s$.  The  coefficient $K_{i,s-i}$  quantifies the coagulation rate and the factor $1/2$ excludes double counting. Similarly, the second term on the r.h.s. of (\ref{SM}) describes the aggregation of particles of size $s$ with all particles of size $i$ with the rates $K_{i,s}$, which leads to the decrease of the concentration $n_s(t)$. Hence, the rates $K_{ij}$ describe binary inter-particle reactions symbolically written as  
$$
[i] + [j] \xrightarrow[]{K_{i,j}} [i+j].  
$$
The dependence of $K_{i,j}$ on $i$ and $j$ is determined  by the physical nature of the process  \cite{Leyvraz, KrapivskyBook, galkin2001smoluchowski, brilliantov2010kinetic, ramkrishna2000population}. 

The solution of the system (\ref{SM}) is challenging, as it is an {\it infinite}  system of equations. Still for the simplest case of a constant kernel $K_{i,j} =1$\footnote{Actually, the constant kernel reads, $K_{i,j}=A$, however the multiplicative constant $A$ may be always chosen as unity by re-scaling of time. The same is true for multiplicative constants in other kernels. } it may be solved for  the mono-disperse initial conditions, $n_s(t = 0) = \delta_{s,1}$\footnote{$\delta_{k,1}$ is the Kronecker symbol.} with the result:
$$
n_i(t)=\frac{4}{(t+2)^2}\left(\frac{t}{t+2}\right)^{j-1}.
$$
This solution was first obtained and analyzed by Smoluchowski \cite{smoluchowski1916drei} and is currently known as the "theory of a rapid coagulation". Later on, a couple of other exact solutions have been obtained for the sum, $K_{i,j}=i+j$, and product, $K_{i,j}=i\cdot j$, kernels \cite{Leyvraz, KrapivskyBook}. Generally, an exact solution is available for  the kernels of the form \cite{Leyvraz, KrapivskyBook}:
$$
K_{i,j}= 1+A(i+j)+ B (i \cdot  j), 
$$
where $A$ and $B$ are constants. 

When the kinetic rates $K_{i,j}$ are homogeneous functions of their arguments $i$ and $j$, there exist also scaling solutions providing the large-time asymptotic $t\to \infty$ for the densities $n_s$.  The most prominent examples of such kernels are  the generalized product kernel, $K_{i,j}=(i\cdot j)^{\mu}$ and  the ballistic kernel \cite{Leyvraz, KrapivskyBook}:
\begin{equation}
 K_{i,j} = (i^{1/3} + j^{1/3})^2 \sqrt{\frac{1}{i} + \frac{1}{j}}.
\end{equation}
Generally, however, Smoluchowski equations may be solved only numerically. 

For many practical applications, a coarse-grained description suffices. In this case the discrete variable $i$ indicating the size of an agglomerate may be treated as continuous and Smoluchowski equations acquires the integro-differential form \cite{muller1928allgemeinen, aloyan1997transport, galkin2001smoluchowski, melzak1957scalar, melzak1957scalarPt2}: 
\begin{equation}
\frac{\partial n(v,t)}{\partial t} = \frac{1}{2} \int_{0}^{v} K(v-u , u) n(v-u,t) n(u,t) du - n(v, t) \int_{0}^{\infty} K(v, u) n(u, t) du.
\end{equation}

Due to the importance of Smoluchowski equations for   industrial applications, such as synthesis of nanoparticles \cite{matveev2020oscillating, boje2017detailed, sabelfeld2018hybrid}, soot formation, pharmacy processes \cite{chaudhury2014computationally}, etc. and for natural phenomena, e.g. \cite{voloshchuk1975coagulation, Falkovich2002, PNAS} a number of  numerical approaches has  been elaborated  \cite{garcia1987monte, liffman1992direct, kruis2000direct, eibeck2000efficient, lee2000validity, MonteCarlo2003, matveev_tensor_2016}. Despite a rapid development of efficient finite-difference and other deterministic methods \cite{stadnichuk_smoluchowski_2015, matveev_fast_2015, osinsky2020low}, the major methodology for solving the aggregation equations remains the stochastic one.  A  particle in cell and Direct Simulation Monte Carlo (DSMC) methods \cite{boje2017detailed, brilliantov2018increasing} may be mentioned as the most prominent examples; they have been proved to be efficient and reliable tools to model aggregation kinetics \cite{ kruis2000direct, sabelfeld2018hybrid, boje2017detailed}.

Early works of Smoluchowski gave rise to the development of  a new field of research associated with agglomeration and related phenomena. The corresponding theoretical and  numerical approaches have been elaborated to describe aggregation with source/sink terms \cite{ZAGAYNOV2001983,zagaynov2007periodic, ball_collective_2012}, multi-particle aggregation \cite{krapivsky1994diffusion, matveev2019numerical}, spontaneous (unary) fragmentation \cite{KrapivskyBook,Bodrova2019}, exchange-driven aggregation \cite{pego2020temporal, niethammer2021oscillations}, etc. One of the most recent extensions of the pure aggregation model was the aggregation-fragmentation model, where the agglomerates undergo collisional (binary) fragmentation  \cite{PNAS,KOB2017,Bodrova2019,stadnichuk_smoluchowski_2015,matveev2017oscillations}. The aggregation-fragmentation kinetics may be  symbolically written for the collisional decomposition into monomers (shattering)  as: 
\begin{align*}
&  [i] + [j] \xrightarrow[]{K_{i,j}} [i+j], & ~~\texttt{aggregation}\\
& [i] + [j] \xrightarrow[]{\lambda \cdot K_{i,j}} \underbrace{[1] + \ldots [1] }_{i+j} & ~~\texttt{shattering},
\end{align*}
where the scalar parameter $\lambda$ quantifies the shattering rate. It may be also shown that the shattering model is generic. Qualitatively, similar results are obtained for more general fragmentation models provided that monomers are predominant in debris size distribution \cite{PNAS}. The corresponding rate equations read, 
\begin{eqnarray}\label{eq:agg_frag_model}
\left\lbrace \begin{matrix}
\dfrac{dn_1}{dt} = & - n_1 \sum\limits_{i=1}^{\infty} K_{1,i} n_i + \lambda \sum\limits_{i = 2}^{\infty} \sum\limits_{j=2}^{\infty} (i+j) K_{i,j} n_i~ n_j + \lambda n_1 \sum\limits_{j=2}^{\infty} j K_{1,j} n_j & \\
 \\
\dfrac{dn_k}{dt} = & \frac{1}{2} \sum\limits_{i=1}^{k-1} K_{i,k-i} n_i n_{k-i} - (1 + \lambda) n_k \sum\limits_{i=1}^{\infty} K_{k,i} n_i, \qquad k = \overline{2,\infty}. \\
\end{matrix} \right. 
\end{eqnarray}
The aggregation-fragmentation equations (\ref{eq:agg_frag_model}) have been  intensively investigated,  analytically and numerically, mainly in the context of particles size distribution in Saturn's rings \cite{stadnichuk_smoluchowski_2015,PNAS}. In these studies, new and efficient deterministic algorithms have been applied; allowing the handling of up to a few hundred thousand equations. With the use of various iterative methods a steady-state and quasi-steady state  size distribution of particles were found \cite{PNAS,stadnichuk_smoluchowski_2015,timokhin2019newton}. The most surprising  was, however,  the detection of never-ending oscillations of particle densities \cite{matveev2017oscillations}. They occur for some range of parameters and may possibly explain the periodic  formation and decay of clumps in the F Ring of Saturn \cite{matveev2017oscillations,French}. The observation of steady oscillations in  closed systems with aggregation and fragmentation was rather surprising, as it was expected that such systems could relax only to a steady or quasi-steady state. Although the numerical evidence of the steady oscillations was quite convincing \cite{PNAS}, an analytical proof  of their existence is still  lacking. The recent efforts to prove theoretically the appearance of the oscillating regime look very promising  \cite{pego2020temporal,niethammer2021oscillations,budzinskiy2020hopf}, yet the problem remains unsolved. 

To prove/disprove the emerging oscillations in aggregation-fragmentation processes is important not only for the fundamental understanding of these processes, but also for their practical applications. Hence it is worth to confirm the existence of density oscillations by an alternative numerical approach, e.g. by stochastic methods. Moreover, stochastic methods would provide an even more solid justification for this phenomenon due to unavoidable presence of noise. Indeed, the deterministic equations (\ref{eq:agg_frag_model})  is essentially an idealized model of a real process. The aggregation-fragmentation kinetics is inevitably accompanied by fluctuations (noise), owing to the discrete nature of microscopic events and a stochastic environment. 

In spite of  the importance of the addressed problem, the stochastic methods have not been widely  applied  to this class of systems yet. To the best of our knowledge, an application of the direct simulation Monte Carlo (DSMC) to aggregation with fragmentation has been reported  in Ref. \cite{sabelfeld2018hybrid} only. Hence, the motivation of the present study is twofold:  Firstly, we extend the existing DSMC methods for aggregation processes with the  collisional fragmentation --  here we probe and adopt two  different DSMC methods \cite{garcia1987monte, kruis2000direct}. Secondly, we apply the extended methods to justify the existence of the  never-ending density oscillations. We confirm that these oscillations are stable and are not sensitive to fluctuations stemming  from the stochastic  nature of DSMC. 


\section{Pure aggregation: two methods revisited}

Applications of DSMC to study aggregation processes has been initiated by Gillespie \cite{gillespie1975exact} who constructed an effective stochastic method for the solution of continuous coagulation equations. During the last decades, the methods have been improved, dramatically increasing the computational speed and accuracy \cite{garcia1987monte, liffman1992direct, matsoukas1997monte, kruis2000direct, eibeck2000efficient, MonteCarlo2003}. Currently, one may classify these methods into several major groups: constant-N Monte Carlo, constant-V Monte Carlo, inverse Monte Carlo, and others. For instance, in the classical case of  binary aggregation, a celebrated  acceptance-rejection method picks randomly  a  pair of particles up and accepts its aggregation with a probability proportional to the aggregation rate \cite{garcia1987monte}. The pseudocode for one step of the acceptance-rejection method is presented in Algorithm \ref{alg:ARmethod}, where $s(i)$ denotes the size of $i$-th particle in the particle array. 
\begin{algorithm}
\caption{Acceptance-rejection method}\label{alg:ARmethod}
\small
\begin{algorithmic}[1]
    \Procedure{ARAggregationStep}{$particles$}
    \State Generate random number $r \in (0,1)$
    \State Choose randomly a pair $(i,j)$
    \State Advance time: $t := t + \tau$,
    \If{$r < \frac{K_{s(i),s(j)}}{K_{\max}}$} 
      \State Add particle of size $s(i) + s(j)$ into array of particles (and delete two old), \label{alg-line:ar-add}
      \State Update $K_{\max}$ and $\tau$ if $s(i) + s(j)$ is larger than the previous largest particle.
    \EndIf
    \EndProcedure
\end{algorithmic}
\end{algorithm}
The computational cost of one step depends on the average aggregation acceptance probability $P_{acc} = \frac{\left\langle K_{s(i), s(j)} \right\rangle}{K_{\max}}$ and the cost of computing $K_{\max}$. If one possesses a complete information about the kernel (symmetry, homogeneity degree, etc.), it is easy to reduce the computational cost for $K_{\max}$ dramatically. Note that $K_{max}$ can be defined as any number satisfying the relation, 
\[
  K_{\max} \geq \mathop {\max }\limits_{i,j} K_{s(i), s(j)}, 
\]
where $s(i)$ and $s(j)$ are currently available particle sizes. It can be limited if the constraint of the  maximal particle size in the system is imposed, and is to be updated only if the size of an emerging aggregate exceeds the current value of $K_{max}$.  If we neglect the cost of $K_{\max}$ computation, then the one-step cost is bounded by the average number of trials before the acceptance, which may be estimated as $1/P_{acc}$. Hence, the cost of a single step is  $O \left( P_{acc}^{-1} \right)$.

An important part of the above algorithm, as well as of  any DSMC method, is the  evaluation rule for  the time-shift $\tau$; it is discussed in detail elsewhere, see e.g. \cite{garcia1987monte, gillespie1975exact,eibeck2000efficient}. For a qualitative understanding of the whole concept, it is worth however,  to sketch the main idea of the rule. 

Let  $n(t) = \sum\limits_{i=1}^{\infty} n_i(t)$ be the total density of particles and $N$ -- their  total number at time $t$, $N=nV$, where $V$ is the system volume. These quantities are related to the initial values at $t=0$ as 
$$
\hat{n}(N) = \frac{N}{N(t = 0)} \,\, n(t=0). 
$$
Consider a small time interval $\Delta t$. According to the definition of the reaction rates, the number of collisions (reactions) between particles of size $i$ and $j$ in the volume $V$ reads, 
$K_{i,j}n_in_j\Delta t V$. The  total number of collisions (reactions) between all particles may be written then as   $(1/2)\sum_{i=1}^{\infty} \sum_{j=1}^{\infty} K_{i,j} n_in_j\Delta t V$. To obtain the collision frequency $\tau_{coll}^{-1}$ one needs to divide the number of collisions by the respective time interval $\Delta t$. Then we obtain:
\begin{eqnarray}
\footnotesize
\label{eq:time_step_assign}
\frac{1}{\tau_{coll}} &=& \frac12 \sum_{i=1}^{\infty} \sum_{j=1}^{\infty} K_{i,j} n_in_j V =
\frac12 \sum_{i=1}^{\infty} \sum_{j=1}^{\infty} \frac{K_{i,j}}{V} N_i N_j 
= \frac{1}{2V} N(N-1) \left\langle K_{s(i), s(j)} \right\rangle \nonumber \\
&=&\frac{K_{\max} ~P_{acc} \cdot N (N - 1)}{2 V} = \frac{K_{\max} ~P_{acc} \cdot \hat n(N) \cdot (N - 1)}{2}, 
\end{eqnarray}
where $N_{i/j} =n_{i/j}V$ and $\left\langle K_{s(i), s(j)} \right\rangle$ is the average reaction rate, which is related to the average acceptance probability as $P_{acc} = \left\langle K_{s(i), s(j)} \right\rangle / K_{\max}$.  Hence  the average time increment $\tau$  may be found as
\[
  \tau = \tau_{coll} \cdot P_{acc}.
\]
Correspondingly, the time-shift for each  trial reads,  
\begin{equation}
\small
\label{eq:time_step_equation}
\tau = \frac{2}{ \hat{n}(N) \cdot (N - 1) \cdot K_{\max}}.
\end{equation}
In our simulations we start with the unit density of particles $n(t = 0) = 1$, which implies that all particles are monomers. 
    
The efficiency of the DSMC can be significantly improved by changing the sampling steps. Instead of sampling a new pair of particles at each time step, one can split this procedure into two separate stages as illustrated in Algorithm \ref{alg:FDSMC}. Firstly, an initial particle is to be chosen for the aggregation event, then a second particle is sampled among the rest $N-1$ particles. This trick allows us to decrease a number of useless rejection events and is known as Fast DSMC (FDSMC).

\begin{algorithm}
\caption{Fast DSMC method}\label{alg:FDSMC}
\small
\begin{algorithmic}[1]
\Procedure{FDSMCAggregationStep}{$particles$}
\State Compute collision frequencies for each particle $S_i = \sum\limits_{j\neq i} K_{s(i), s(j)}$.
\State Choose random $r \in \left(0, \sum\limits_{k=1}^{N}S_i\right)$.
\State Find $i$ s.t. $\sum\limits_{k = 1}^{i - 1}S_k \leq r \leq \sum\limits_{k = 1}^{i}S_k$  \Comment{first particle selection}
\State Find $j$ s.t. $\sum\limits_{k = 1}^{j - 1}K_{s(i), s(k)} \leq  r - \sum\limits_{k = 1}^{i - 1}S_k \leq \sum\limits_{k = 1}^{j}K_{s(i), s(k)}$ \Comment{second particle selection}
\State Create particle  $s(i) + s(j)$ and delete $s(i)$ and $s(j)$. \label{alg-line:fdsmc-add}
\State Advance time $t := t + \tau$.
\EndProcedure
\end{algorithmic}
\end{algorithm}
This trick requires to compute, store and efficiently update the partial sums $S(i)$. 
Additional acceleration of such types of methods can be achieved by using majorants of the kernel coefficients to approximate the evaluations of the sums (see \cite{eibeck2000efficient}). One can also apply grouping of particles of the fixed size into the ``buckets'' and achieve extra speedup; this allows to deal with billions of particles. 

Noteworthy, DSMC can be  rather efficiently deployed on parallel computers. The according  theoretical analysis can be found in Ref. \cite{MonteCarlo2003} and the discussion of the experimental implementation in Ref. \cite{xu2014fast} for CPU-based clusters and in Refs. \cite{wei2013gpu, xu2015accelerating} for modern GPUs.

Regardless of the chosen algorithm, a user has to decide what to do when only a few particles are left in the system. One possible and very natural choice is to clone all sampled particles from time to time without a significant loss of accuracy \cite{liffman1992direct}. 

Certainly, cloning of particles can be done in many ways \cite{matsoukas1997monte,BFP2018}. The most straightforward way is to increase the volume of the system --  this keeps the number of particles approximately constant; the simulation accuracy is also kept  at the target level. For instance, one can duplicate the amount of the existing particles once their number drops down below half of the initial number. This mimics the doubling of the volume of the simulation cell, conserving the total mass density.

\section{Application to collisional shattering}

The modification of both methods for the case of (binary) collisional shattering may be done by taking into account that the particle array updates are decoupled from the pair selection procedure. Since the shattering rate $\lambda$ does not depend on the particle densities and their size, we can choose the type of interaction -- aggregation or fragmentation, afterwards, after the pair has been picked up. 

Algorithms \ref{alg:ARmethodMOD} and \ref{alg:FDSMC_frag} 
illustrate the extension of the acceptance-rejection and FDSMC algorithms by incorporating binary shattering process. The decision on the type of process (aggregation or fragmentation) chosen is included as a separate procedure \verb$PairInteraction$. The procedure can be used as a drop-in replacement for updates in lines \ref{alg-line:ar-add} in both Algorithms \ref{alg:ARmethod} and \ref{alg:FDSMC}. These  additional pairwise interactions may be directly added to the procedure in an obvious manner.

\begin{algorithm}
\caption{Acceptance-rejection method for binary scattering}\label{alg:ARmethodMOD}
\small
\begin{algorithmic}[1]
    \Procedure{ARAggregationFragmentationStep}{$particles$}
    \State Generate random number $r \in (0,1)$
    \State Choose random pair $(i,j)$
    \State Advance time: $t := t + \tau$,
    \If{$r < \frac{K_{s(i),s(j)}}{K_{\max}}$} 
      \State PairInteraction($i$, $j$) \label{alg-line:ar-interaction}
    \EndIf
    \EndProcedure
    
    \Procedure{PairInteraction}{$i$, $j$}
    \State Generate random number $r \in (0,1)$ 
    \If {$r < \frac{\lambda}{1 + \lambda}$} \Comment{fragmentation}
      \State Add $s(i) + s(j)$ monomers into array of particles.
    \Else \Comment{aggregation}
      \State Add particle of size $s(i) + s(j)$ into array of particles.
      \State Update $K_{\max}$ if $s(i) + s(j)$ is larger than the previous largest particle.
    \EndIf
    \State Recompute $\tau$.
    \State Delete particles $i$ and $j$.
    \EndProcedure
\end{algorithmic}
\end{algorithm}

The other significant change in this modification is the updated rule for time-shift calculation, which now reads, 
\begin{equation}
\label{eq:AR_fragmentation}
\tau = \frac{2}{\hat{n}(N) (N-1) K_{\max}} \cdot \frac{1}{1 + \lambda}.
\end{equation}
The same factor $\frac{1}{1 + \lambda}$ has to be used to update the time-shifts in the FDSMC algorithm.

\begin{algorithm}
\caption{Fast DSMC method for binary shattering}\label{alg:FDSMC_frag}
\small
\begin{algorithmic}[1]
\Procedure{FDSMCAggregationFragmentationStep}{$particles$}
\State Compute collision frequencies for each particle $S_i = \sum\limits_{j\neq i} K_{s(i), s(j)}$.
\State Choose random $r \in \left(0, \sum\limits_{k=1}^{N}S_i\right)$.
\State Find $i$ s.t. $\sum\limits_{k = 1}^{i - 1}S_k \leq r \leq \sum\limits_{k = 1}^{i}S_k$  \Comment{first particle selection}
\State Find $j$ s.t. $\sum\limits_{k = 1}^{j - 1}K_{s(i), s(k)} \leq  r - \sum\limits_{k = 1}^{i - 1}S_k \leq \sum\limits_{k = 1}^{j}K_{s(i), s(k)}$ \Comment{second particle selection}
\State PairInteraction($i$, $j$). \Comment{reused from Algorithm \ref{alg:ARmethodMOD}} \label{alg-line:fdsmc-interaction}
\State Advance time $t := t + \tau$.
\EndProcedure
\end{algorithmic}
\end{algorithm}

While Algorithms \ref{alg:ARmethodMOD} and \ref{alg:FDSMC_frag}  rely on two independent random numbers for the pair selection and then for the event selection, it is possible to use only one random number,  by generating $r$ uniformly distributed in $(0, 1 + \lambda)$. The range $(0, 1)$ then corresponds to the aggregation event and the  range $(1, 1 + \lambda)$ to the fragmentation event. 

\section{Particle grouping}

The following  feature of the shattering  algorithms hinders their  efficiency -- the need  to store a large number of newly created monomers right after the fragmentation event.

Using the fact that particles of the same size are indistinguishable, we suggest  to group the particles into ``buckets''. Each bucket comprises particles of one size only and keeps their number along with other information, such as total probability $S(i)$. We store the buckets in a static array allocated ahead in time. Each index of the array  directly corresponds to the particle size, allowing for a quick search of the bucket for a correct size during the insertions.
To support the particles' unbounded growth and avoid a huge number of empty buckets in the preallocated array, we  group only particles, whose size does not exceed $B$. If a particle is larger than $B$, it is added to a dynamic vector as a separate particle.   

With this hybrid static-dynamic approach, the addition of the monomers requires just an update of the counter in the bucket and a single traversal over all buckets and particles in the dynamic array; this allows  to update all interaction probabilities.

Tuning the parameter $B$, one achieves a trade-off between a cost of bucket array traversal and a cost of storing particles as separate units. In the extreme case of $B \to \infty$, most of the buckets will be empty and any particle insertion/removal will incur a cost proportional to the size of the biggest particle. Another extreme case $B = 0$ corresponds to the lack of grouping, and any particle update will incur a cost proportional to the number of particles, as stated in the original algorithm.  The choice of $B \sim 1000$ guarantees a good performance for most of the practical applications. 

\section{Numerical experiments}
\subsection{Performance and accuracy testing on the basic problems}

The very first test measures the accuracy of the aggregation-shattering Monte Carlo algorithms for the case of the constant kernel,  
$$
K_{i,j} = 1
$$
and monodisperse initial conditions $n_s(t = 0) = \delta_{s,1}$; an analytical solution for this case is available, see e.g. \cite{PNAS}. In Table \ref{table:eucl_constant}, we present the results for the accuracy of the method,  quantified by the Euclidean norm of the error of the final particle size distribution  (recall that the exact solution is known). The results indicate a good accuracy of the method and convergence to the analytical solution. No significant difference between the acceptance-rejection and FDSMC schemes is observed. Noteworthy, the error converges as $O(N^{-1/2})$. 

\begin{figure}[t]
\centering
\includegraphics[width=1.0\columnwidth]{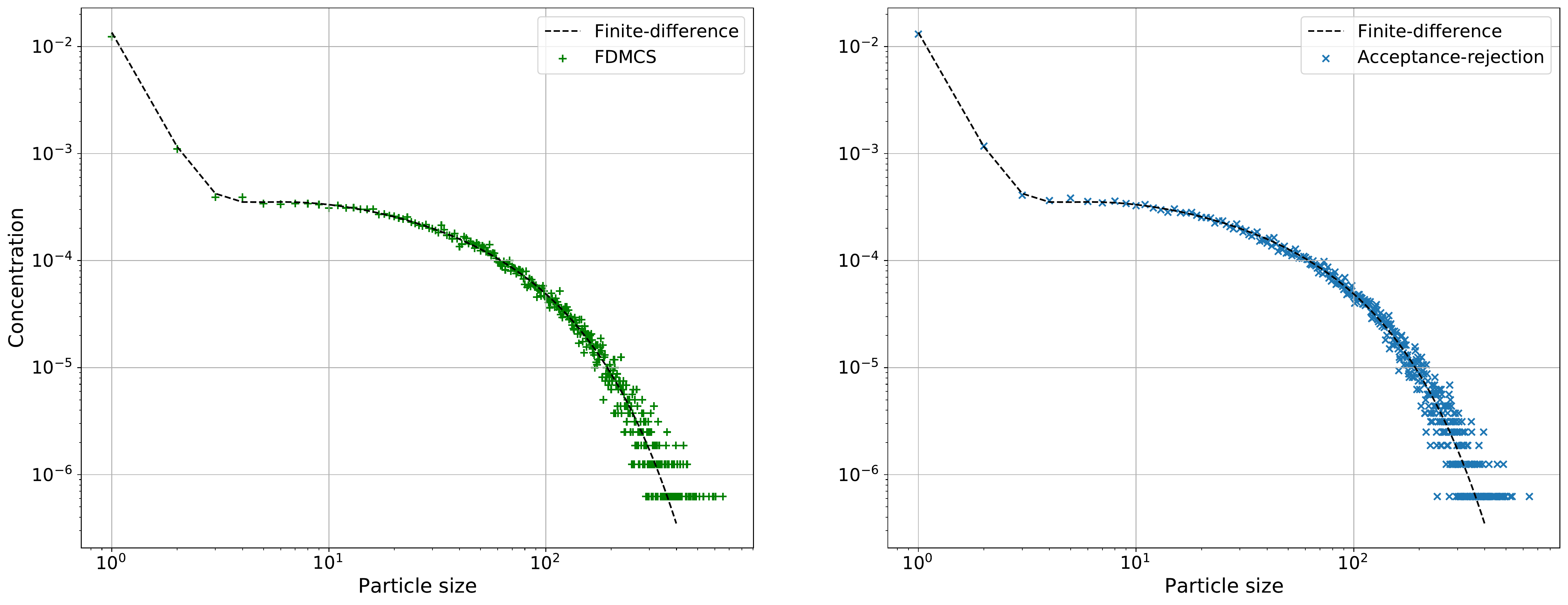}
\caption{Particle size distribution at $t=10$ for the ballistic kernel obtained by MC algorithms for the initial number of particles $N=10^5$ (dots) and by the finite-difference method (lines) with the step size $\Delta t = 10^{-3}$ and $5 000$ equations. The shattering rate  $\lambda = 0.01$ and mono-disperse initial conditions are used. Right panel: Fast DSMC. Left panel: The acceptance-rejection method. }
\label{fig:ballistic_separate}
\end{figure}

\begin{table}[]
\begin{center}
\small
\begin{tabular}{l|ll|ll|ll}
           & \multicolumn{2}{c|}{$t=1$} & \multicolumn{2}{c|}{$t=10$} & \multicolumn{2}{c}{$t=100$} \\ \cline{2-7} 
           & A-R     & FDSMC   & A-R      & FDSMC    & A-R     & FDSMC \\ \hline
$N=10^3$   & 0.012   &  0.010  &  0.008   &  0.019   & 0.017   & 0.034 \\
$N=10^4$   & 0.005   & 0.003   &  0.009   &  0.004   & 0.009   &  0.008      \\
$N=10^5$   & 0.002   & 0.0015   & 0.0007   &   0.001  & 0.0023 &  0.0024 \\
$N=10^6$   & 0.0004  &  0.0003   &  0.0005 &  0.0005  & 0.0012 &  0.0005\\
\hline 
\end{tabular}
\end{center}
\caption{The Euclidean norm of the error of the solution obtained by the acceptance-rejection (A-R) and FDMSC algorithms with the constant kernel, monodisperse initial conditions and $\lambda = 0.1$.}
\label{table:eucl_constant}
\end{table}

\begin{table}[]
\begin{center}
\scriptsize
\begin{tabular}{l|l|lll|lll|lll}
Method & Parame- & \multicolumn{3}{c|}{$t=1$} & \multicolumn{3}{c|}{$t=10$} & \multicolumn{3}{c}{$t=100$} \\ \cline{3-11}
    &  ters       & $N(t)$ & $M_2(t)$ &  $M_3(t)$ & $N(t)$ & $M_2(t)$ &  $M_3(t)$ & $N(t)$ & $M_2(t)$ &  $M_3(t)$ \\ \hline

    {\tiny A-R} & \multirow{2}{*}{{\tiny $N=10^3$}} & 3\% & 5\% & 16\% & 26\% & 8\% & 21\% & 42\% & 28\% & 49\% \\
    {\tiny FDSMC} & & 16.7\% & 0.7\% & 1.3\% & 65\% & 0.4\% & 2.9\% & 21\% & 5\% & 2.1\% \\ \hline

    {\tiny A-R} & \multirow{2}{*}{{\tiny $N=10^4$}} & 0.9\% & 2\% & 5\% & 10\% & 2\% & 7\% & 13\% & 9\% & 14\% \\
    {\tiny FDSMC} & & 2.2\% & 0.57\% & 3.14\% & 78.5\% & 0.13\% & 1.34\% & 6.6\% & 16.6\% & 5.8\% \\ \hline

    {\tiny A-R} & \multirow{2}{*}{{\tiny $N=10^5$}} & 0.2\% & 0.6\% & 2\% & 4\% & 0.9\% & 2\% & 5\% & 3\% & 6\% \\
    {\tiny FDSMC} & & 0.8\% & 0.32\% & 0.65\% & 9.6\% & 0.9\% & 2.9\% & 5.9\% & 1.5\% & 1.3\% \\ \hline

    {\tiny A-R} & \multirow{2}{*}{{\tiny $N=10^6$}} & 0.08\% & 0.2\% & 0.6\% & 0.6\% & 0.2\% & 0.6\% & 2\% & 0.7\% & 2\% \\
    {\tiny FDSMC} & & 0.38\% & 0.02\% & 0.08\% & 0.4\% & 0.2\% & 0.7\% & 1.1\% & 1.9\% & 2.3\% \\ \hline
    
\multirow{2}{*}{{\tiny FD}}

    & {\tiny $\Delta t = 0.05$} & 3.8\% & 4\% & 8.1\% & 0.17\% & 0.46\% & 0.92\% & 0.1 \% & 0.15 \% & 0.66 \% \\
    & {\tiny $\Delta t = 0.005$} & 0.33\% & 0.4\% & 0.81\% & 0.014\% & 0.05\% & 0.09\% &  0.0016 \% & 0.0017 \% & 0.0033 \% \\ \hline
    
\end{tabular}
\end{center}
\caption{
Comparison of the errors of the simulation results for  zeroth $M_0(t) = \sum\limits_{s=1}^{\infty} n_s(t)$, second $M_2(t) = \sum\limits_{s=1}^{\infty} s^2 n_s(t)$ and third $M_3(t) = \sum\limits_{s=1}^{\infty} s^3 n_s(t)$ moments of the particle size distribution for the ballistic kernel and $\lambda = 0.01$. The average relative error (standard deviation) is shown for both MC algorithms and the  finite-difference scheme. Different number of particles (for the MC algorithms) and different time steps (for the finite-difference scheme) are used.  }
\label{tab:ballistic_rel_convergence}
\end{table}

Next we report the benchmarks for the case of the ballistic kernel,
$$
K_{i,j} = (i^{1/3} + j^{1/3})^2 \sqrt{\frac{1}{i} + \frac{1}{j}}, 
$$
with the same initial conditions. For this problem, the analytical solution is not known and the kernel itself requires more ``heavy'' computations.  In Table \ref{tab:ballistic_rel_convergence} and \ref{subfig:excintction-a} we demonstrate the convergence of the results for the total density and two higher-order moments of the particle size distribution:  $$n(t) = \sum\limits_{s=1}^{\infty} n_s(t) \qquad M_2(t) = \sum\limits_{s=1}^{\infty} s^2 n_s(t) \qquad
M_3(t) = \sum\limits_{s=1}^{\infty} s^3 n_s(t). 
$$
The data in Table \ref{tab:ballistic_rel_convergence} clearly illustrates that the accuracy of the finite-difference approach is higher. 

In Tables \ref{table:runtime_constant} and \ref{table:runtime_ballistic} we present the CPU times of simulations by the DSMC algorithms for the constant and ballistic kernels. It is compared with the according simulation time for the fast finite-difference approach \cite{matveev_fast_2015}. One can see that although the  DSMC approach is less accurate, it provides faster simulations (see Table \ref{table:runtime_ballistic})  ensuring the conservation of mass. 

\begin{table}[]
\begin{center}
\small
\begin{tabular}{l|ll|ll|ll}
           & \multicolumn{2}{c|}{$t=1$} & \multicolumn{2}{c|}{$t=10$} & \multicolumn{2}{c}{$t=100$} \\ \cline{2-7} 
           & A-R    &  FDSMC  &  A-R    &  FDSMC  & A-R    &  FDSMC \\ \hline
$N=10^3$   & 0.0001   &  0.007  &  0.0001   &  0.07   & 0.0001   &  0.64 \\
$N=10^4$   & 0.0001   &  0.005  &  0.0006   &  0.09   & 0.0025   &  0.82      \\
$N=10^5$   & 0.001   &  0.046  &  0.0059  &  0.13   & 0.025 &  1.97 \\
$N=10^6$   & 0.011  &  0.09   &  0.046 &  0.71   & 0.26 &  15.13 \\ \hline
Finite-difference & \multicolumn{2}{c|}{$2.49$} & \multicolumn{2}{c|}{$24.66$} & \multicolumn{2}{c}{$246.21$} \\
\hline
\end{tabular}
\end{center}
\caption{Running time, in seconds, for  the acceptance-rejection (A-R), Fast Direct Monte Carlo Simulation (FDSMC) and Finite-Difference ($32~758$ equations) algorithms for the constant kernel, monodisperse initial conditions and $\lambda = 0.1$.}
\label{table:runtime_constant}
\end{table}

\begin{table}[]
\begin{center}
\small
\begin{tabular}{l|ll|ll|ll}
           & \multicolumn{2}{c|}{$t=1$} & \multicolumn{2}{c|}{$t=10$} & \multicolumn{2}{c}{$t=100$} \\ \cline{2-7} 
           & A-R     & FDSMC   &  A-R     & FDSMC    & A-R    &  FDSMC \\ \hline
$N=10^3$   & 0.0001   &  0.018  &  0.0005   &  0.407   & 0.071   &  245.4 \\
$N=10^4$   & 0.0011    &  0.15   &  0.0064   &  3.4     & 1.11   &  439.0       \\
$N=10^5$   & 0.013     &  1.67   &  0.076    &  34.82   & 14.15      &  5204 \\
$N=10^6$   & 0.17   &  18.08  &   0.96      &  406.8   & 87.89      &  40335   \\ \hline
Finite-difference & \multicolumn{2}{c|}{$9.51$} & \multicolumn{2}{c|}{$94.57$} & \multicolumn{2}{c}{$988.48$}\\
\hline
\end{tabular}
\end{center}
\caption{Running time, in seconds, for the acceptance-rejection (A-R), Fast Direct Monte Carlo Simulation (FDSMC) and Finite-Difference algorithm ($5~000$ equations) for the  ballistic kernel, monodisperse initial conditions and $\lambda = 0.01$.}
\label{table:runtime_ballistic}
\end{table}


\subsection{Oscillating regimes in the aggregation-fragmentation kinetics}

\begin{figure}[ht]
    \centering
    \begin{subfigure}{0.48\textwidth}
        \includegraphics[width=\columnwidth]{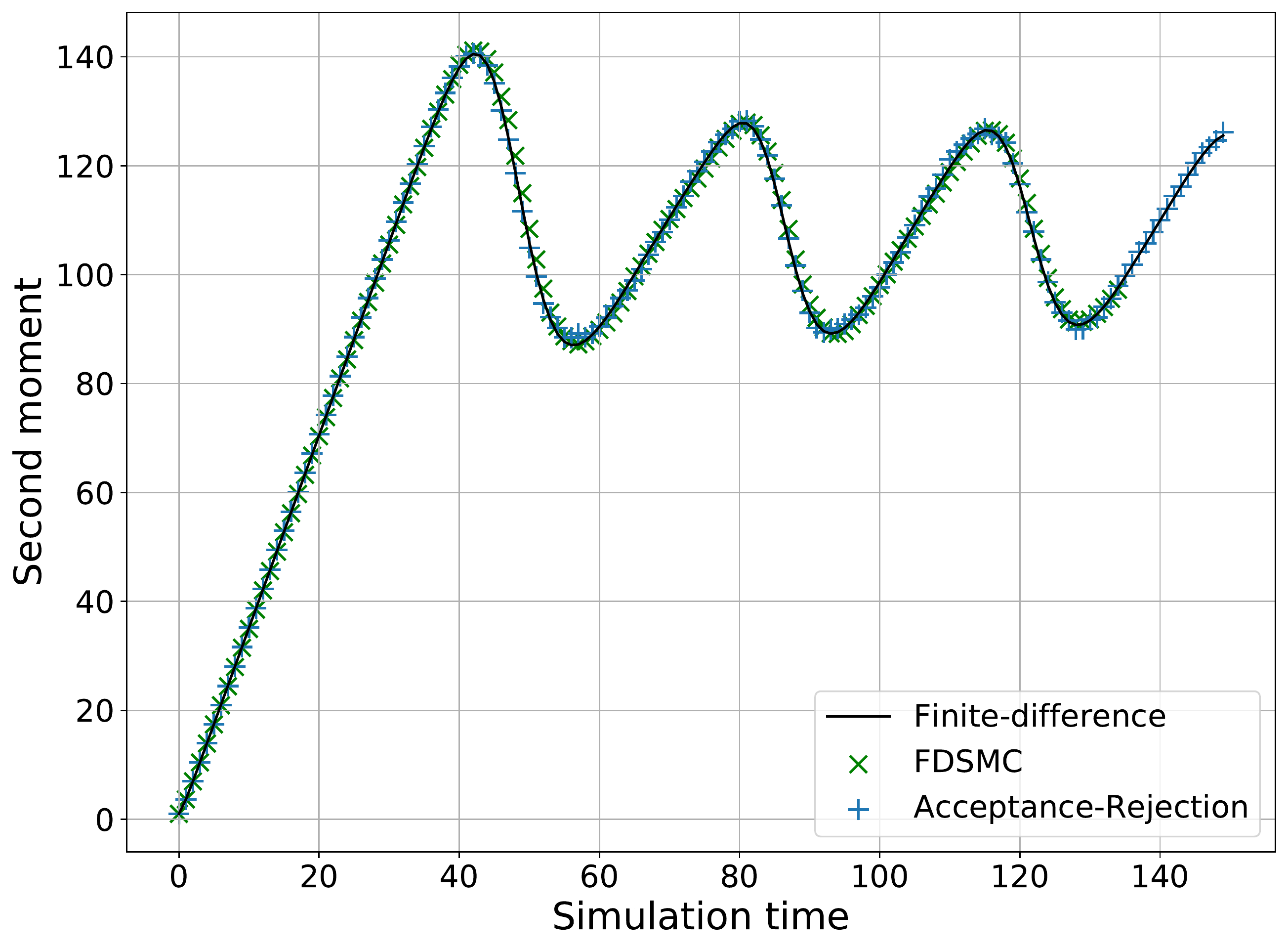}
        \caption{}
        \label{fig:osc_a}
    \end{subfigure}
    \begin{subfigure}{0.48\textwidth}
        \centering
        \includegraphics[width=\columnwidth]{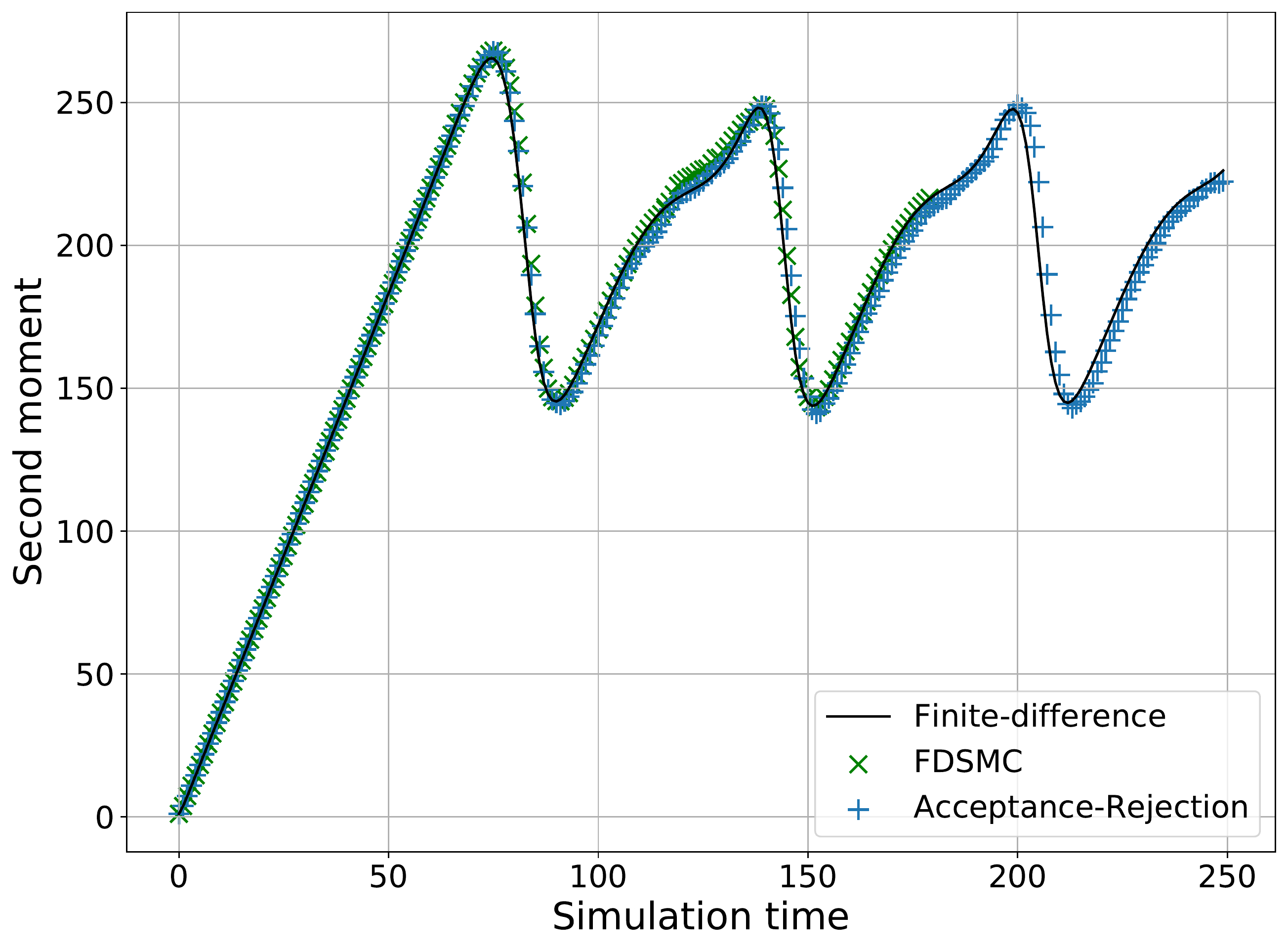}
        \caption{}
        \label{fig:osc_b}
    \end{subfigure}
    \caption{Time-dependence of the second moment of the particle size distribution function, $M_2(t) = \sum\limits_{k=1}^{\infty} k^2 n_k(t)$ for the case of monodisperse initial conditions. (a) $a = 0.95$ and $\lambda = 0.01$, (b) $a = 0.98$ and $\lambda = 0.005$. Dots -- Monte Carlo results, lines -- the results of the finite-difference solution of $8192$ equations. One can see the emergence of steady oscillations. } 
    \label{fig:Oscillations}
\end{figure}
Opposite to the plain relaxation behavior expected for the Smoluchowski-like equations, a surprising oscillation regime has  been recently reported for the aggre\-gation-fragmentation  kinetics. It was observed that the kinetic equations \eqref{eq:agg_frag_model} with the kernels 
\begin{equation*}
    K_{i,j} = \left(\frac{i}{j}\right)^{a} + \left(\frac{j}{i}\right)^{a}
\end{equation*}
demonstrated for closed systems, with a lack of sources and  sinks of particles,  steady oscillations for  $a > 0.5$ and $0 < \lambda < \lambda_{crit}$ \cite{matveev2017oscillations,BrillPRE2018}. Such oscillations have been detected numerically, by an efficient implementations of the finite-difference second-order Runge-Kutta time-integration scheme. Although the qualitative theory of Ref. \cite{BrillPRE2018} supports  the existence of the oscillating regimes, an additional independent verification of this phenomenon is highly desirable. Here we show  that DSMC methods demonstrate the same oscillations of densities, as  in Refs. \cite{matveev2017oscillations,BrillPRE2018} and the according numerical solutions converge, see Figure \ref{fig:Oscillations}.   Both, the acceptance-rejection and FDSMC approach demonstrate the emergent oscillations, thus verifying the results of the previous studies. Therefore, we conclude that the oscillatory solutions are robust and stable with respect to stochastic errors and fluctuations inherent for Monte Carlo methods. 

\begin{figure}
    \centering
    \begin{subfigure}{0.48\textwidth}
    \centering
    \includegraphics[width=\columnwidth]{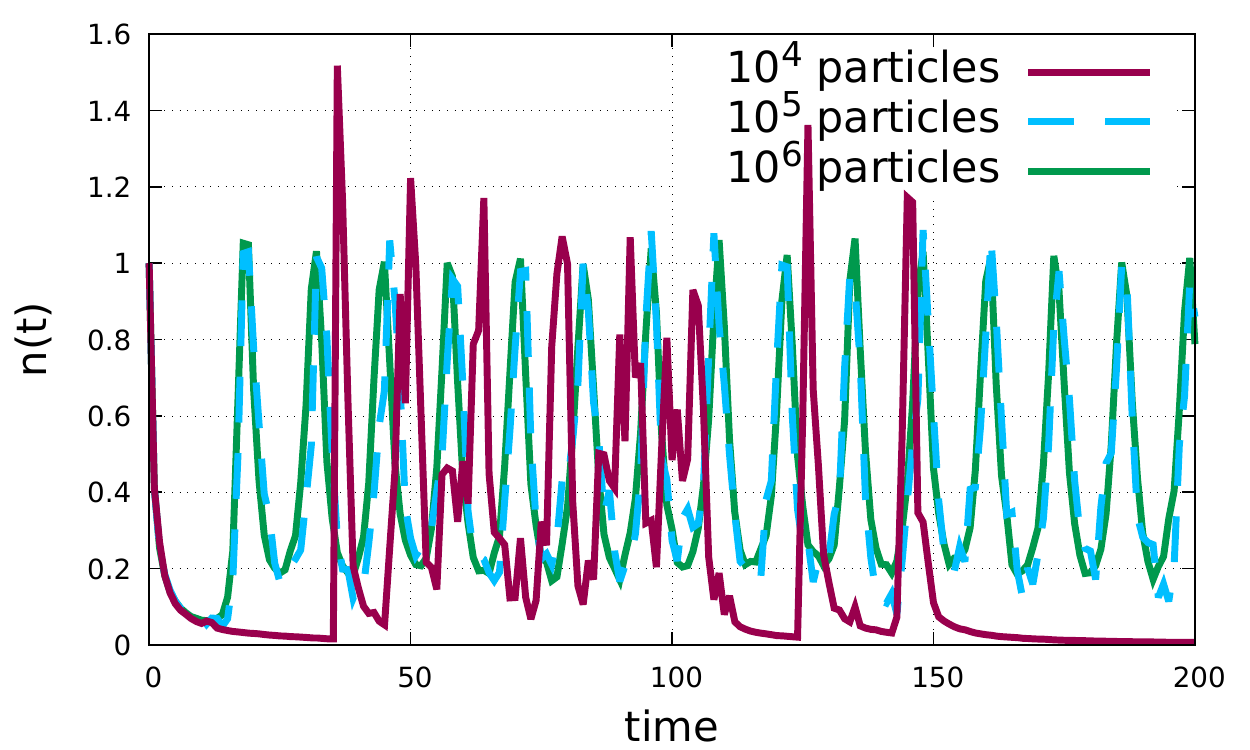}
    \caption{}
    \label{subfig:excintction-a}
    \end{subfigure}
    \begin{subfigure}{0.48\textwidth}
    \centering
    \includegraphics[width=\columnwidth]{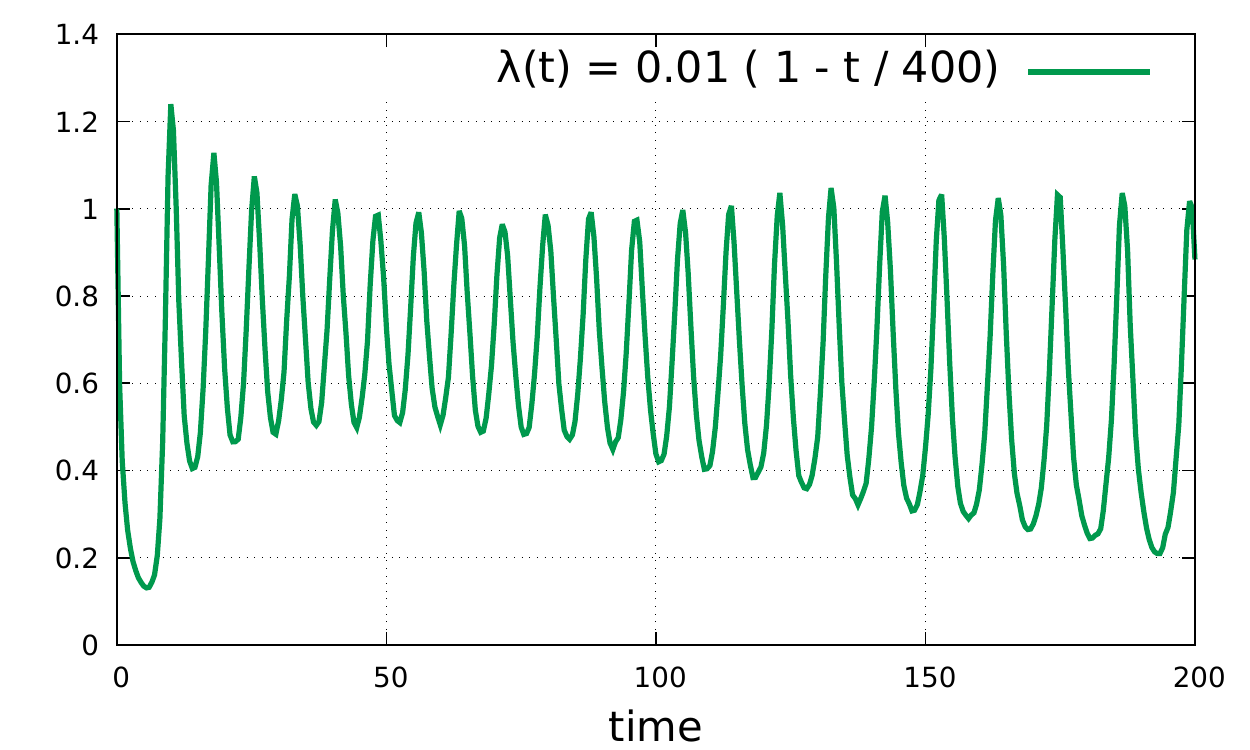}
    \caption{}
        \label{subfig:timdeplambda}
    \end{subfigure}
    \caption{(a) Demonstration of an ``extinction'' effect for the case $a = 0.95$, $\lambda = 0.005$ and polydisperse initial conditions $n_k(t=0)=0.1, k = 1, \ldots 10$, $n_k(t=0) = 0, k > 10$. One may see the convergence of results as number of particles in simulation grows but for $10^4$ (dashed line) we obtain the ``extinction'' effect at a large time. (b) Simulations for time-dependent $\lambda(t) = 0.01\cdot (1-t/400)$, the starting value of $\lambda(t=0) = 0.01$ corresponds to the scenario of damping oscillations which should converge to the fixed point but the decreasing value of shattering causes the solution switch to more pronounced and steady oscillations instead of the steady-state.}
    \label{fig:extinction}
\end{figure}

In our experiments, we also observed an  ``extinction'' effect  (see Fig. \ref{subfig:excintction-a}), when all particles coalesce in a single final cluster; evolution of the systems ceases in this case. The extinction occurs for  relatively long sequences of only-aggregative collisions. The probability of such sequences rapidly decreases with the systems size.  Hence to observe the steady oscillations the system should be large enough. Still for any finite system the extinction effect is always possible, even if its probability is vanishingly  small.

Note that the possibility of this effect illustrates a very important difference between the  stochastic  Monte Carlo approach and the deterministic one  of the finite difference method. The latter approach ignores the discrete nature of the process with inevitable fluctuations, and essentially models an infinitely large system in the thermodynamic limit.  Hence the reported results confirm the existence of steady oscillations in aggregation-fragmentation kinetics in real systems, which are always finite and subjected to noise and fluctuations. 

Another interesting observation comes from experiments with time-dependent shattering parameter $\lambda(t)$. In Fig. \ref{subfig:timdeplambda} for $a = 0.9$ we demonstrate that damping oscillations (e.g. for $a = 0.9$, $\lambda = 0.01$ we know that the solution converges to the steady-state \cite{matveev2017oscillations}) can become steady if $\lambda$ smoothly decreases to the region corresponding to oscillations. In our simulations from Fig. \ref{fig:extinction} we used the linear relaxation
$$\lambda(t) = 0.01\cdot \left(1 - \frac{t}{400}\right)$$
for $t \in [0,200]$ and the same initial conditions as in our previous work \cite{matveev2017oscillations}:
\begin{equation}\notag
    n_k(t=0) = \left\lbrace \begin{matrix} 
    0.1, & k = 1, 2, \ldots 10 \\
    0, & k > 10 \\
    \end{matrix}\right.
\end{equation}
These results allow us to propose a conjecture that oscillations in this class of systems rise through a Hopf bifurcation by analogy with recent results \cite{pego2020temporal, budzinskiy2020hopf} for simpler models addition-shattering and exchange kinetics. We hope to provide more detailed analysis of this phenomena in our future research.

\section{Conclusions}

In the present study we report a straightforward generalization of two well-established Monte Carlo (MC) techniques for aggregation kinetics, supplemented by collisional fragmentation (shattering).  We present a comprehensive  validation of the proposed methodology. The MC results are compared with the known analytical solutions for the constant kernel and with the alternative efficient finite-difference scheme \cite{matveev_fast_2015} for the ballistic kernel. We show that the basic MC methods -- the acceptance-rejection and FDSMC method,  formulated initially for a pure aggregation may be easily modified for a more complex model without losing  the computational efficiency and good convergence. We implement both methods using the grouping of particles and demonstrate that it  significantly outperforms,  in terms of the CPU-time,  the most efficient finite-difference method of Ref. \cite{matveev_fast_2015}.

We also demonstrate that the recent observations of the oscillatory regimes in the aggregation-fragmentation kinetics \cite{matveev2017oscillations,BrillPRE2018} for the generalized Brownian kernel may  be validated by the stochastic methods. Thus, our experiments provide a convincing proof that the steady  oscillations do exist and are stable with respect to fluctuations and noise. 

\section*{Acknowledgements}

Zhores supercomputer of Skolkovo Institute of Science and Technology \cite{zacharov2019zhores} has been used in the present research. S.M. was supported by Moscow Center for Fundamental and Applied Mathematics (the agreement with the Ministry of Education and Science of the Russian Federation No. 075-15-2019-1624). A.O. acknowledges RFBR project No. 20-31-90022 and N.B. -- RFBR project No. 18-
29-19198. 

\bibliography{Smoluchowski}

\begin{thebibliography}{10}

\bibitem{Leyvraz}
F.~Leyvraz.
\newblock Scaling theory and exactly solved models in the kinetics of
  irreversible aggregation.
\newblock {\em Phys. Reports}, 383:95--212, 2003.

\bibitem{KrapivskyBook}
P.~L. Krapivsky, S.~Redner, and E.~Ben-Naim.
\newblock {\em A Kinetic View of Statistical Physics}.
\newblock Cambridge University Press, 2010.

\bibitem{Krapivsky_PRB1999}
P.~L. Krapivsky, J.~F.~F. Mendes, and S.~Redner.
\newblock Influence of island diffusion on sub-monolayer epitaxial growth.
\newblock {\em Physical Review B}, 59:15950, 1999.

\bibitem{SpahnetalEPL:2004}
F.~{Spahn}, N.~{Albers}, M.~{Sremcevic}, and C.~{Thornton}.
\newblock {Kinetic description of coagulation and fragmentation in dilute
  granular particle ensembles}.
\newblock {\em Europhysics Letters}, 67:545--551, 2004.

\bibitem{BrilliantovSpahn2006}
N.~V. Brilliantov and F.~Spahn.
\newblock Dust coagulation in equilibrium molecular gas.
\newblock {\em Mathematics and Computers in Simulation}, 72:93, 2006.

\bibitem{DominikTilens:1997}
Dominik C. and Tielens A.~G. G.
\newblock The physics of dust coagulation and the structure of dust aggregates
  in space.
\newblock {\em Astrophys. J.}, 480:647, 1997.

\bibitem{PNAS}
N.~V. Brilliantov, P.~L. Krapivsky, A.~Bodrova, F.~Spahn, H.~Hayakawa,
  V.~Stadnichuk, and J.~Schmidt.
\newblock Size distribution of particles in {S}aturn's rings from aggregation
  and fragmentation.
\newblock {\em PNAS}, 112(31):9536--9541, 2015.

\bibitem{BFP2018}
N.~Brilliantov, A.~Formella, and T.~Poeschel.
\newblock Increasing temperature of cooling granular gases.
\newblock {\em Nature Communications}, 9:797, 2018.

\bibitem{prions}
T.~Poeschel, N.~V. Brilliantov, and C.~Frommel.
\newblock Kinetics of prion growth.
\newblock {\em Biophys. J.}, 85:3460--3474, 2003.

\bibitem{CoagNetw}
W.~Miura, H.~Takayasu, and M.~Takayasu.
\newblock Effect of coagulation of nodes in an evolving complex network.
\newblock {\em Phys. Rev. Lett}, 108:168701, 2012.

\bibitem{SA08}
K.~Ariga, J.~P. Hill, M.~V. Lee, A.~Vinu, R.~Charvet, and S.~Acharya.
\newblock Challenges and breakthroughs in recent research on self-assembly.
\newblock {\em Sci. Tech. Adv. Mater.}, 9:014109, 2008.

\bibitem{Igor14}
A.~Demortire, A.~Snezhko, M.~V. Sapozhnikov, N.~Becker, T.~Proslier, and I.~S.
  Aranson.
\newblock Self-assembled tunable networks of sticky colloidal particles.
\newblock {\em Nat. Commun.}, 5:3117, 2014.

\bibitem{Erik17}
C.~G. Evans and E.~Winfree.
\newblock Physical principles for dna tile self-assembly.
\newblock {\em Chem. Soc. Rev.}, 46:3808, 2017.

\bibitem{smoluchowski1916drei}
M.~V. Smoluchowski.
\newblock Drei vortrage uber diffusion, {B}rownsche bewegung und koagulation
  von kolloidteilchen.
\newblock {\em Zeitschrift fur Physik}, 17:557--585, 1916.

\bibitem{galkin2001smoluchowski}
V.A. Galkin.
\newblock Smoluchowski equation.
\newblock {\em Fizmatlit, Moscow}, 2001.

\bibitem{brilliantov2010kinetic}
Nikolai~V Brilliantov and Thorsten P{\"o}schel.
\newblock {\em Kinetic theory of granular gases}.
\newblock Oxford University Press, 2010.

\bibitem{ramkrishna2000population}
Doraiswami Ramkrishna.
\newblock {\em Population balances: {Theory and applications to particulate
  systems in engineering}}.
\newblock Elsevier, 2000.

\bibitem{muller1928allgemeinen}
H.~M{\"u}ller.
\newblock Zur allgemeinen theorie ser raschen koagulation.
\newblock {\em Fortschrittsberichte {\"u}ber Kolloide und Polymere},
  27(6):223--250, 1928.

\bibitem{aloyan1997transport}
A.~E. Aloyan, V.~O. Arutyunyan, A.~A. Lushnikov, and V.~A. Zagaynov.
\newblock Transport of coagulating aerosol in the atmosphere.
\newblock {\em Journal of {Aerosol} {Science}}, 28(1):67--85, 1997.

\bibitem{melzak1957scalar}
ZA~Melzak.
\newblock A scalar transport equation.
\newblock {\em Transactions of the American Mathematical Society},
  85(2):547--560, 1957.

\bibitem{melzak1957scalarPt2}
ZA~Melzak et~al.
\newblock A scalar transport equation. {I}{I}.
\newblock {\em The Michigan Mathematical Journal}, 4(3):193--206, 1957.

\bibitem{matveev2020oscillating}
SA~Matveev, AA~Sorokin, AP~Smirnov, and EE~Tyrtyshnikov.
\newblock Oscillating stationary distributions of nanoclusters in an open
  system.
\newblock {\em Mathematical and Computer Modelling of Dynamical Systems}, pages
  1--14, 2020.

\bibitem{boje2017detailed}
Astrid Boje, Jethro Akroyd, Stephen Sutcliffe, John Edwards, and Markus Kraft.
\newblock Detailed population balance modelling of {T}i{O}2 synthesis in an
  industrial reactor.
\newblock {\em Chemical Engineering Science}, 164:219--231, 2017.

\bibitem{sabelfeld2018hybrid}
Karl~K Sabelfeld and Georgy Eremeev.
\newblock A hybrid kinetic-thermodynamic {M}onte {C}arlo model for simulation
  of homogeneous burst nucleation.
\newblock {\em Monte Carlo Methods and Applications}, 24(3):193--202, 2018.

\bibitem{chaudhury2014computationally}
A.~Chaudhury, I.~Oseledets, and R.~Ramachandran.
\newblock A computationally efficient technique for the solution of
  multi-dimensional {P}{B}{M}s of granulation via tensor decomposition.
\newblock {\em Computers \& Chemical Engineering}, 61:234--244, 2014.

\bibitem{voloshchuk1975coagulation}
V.~M. Voloshchuk and Y.~S. Sedunov.
\newblock {\em Coagulation processes in disperse systems}.
\newblock Gidrometeoizdat, Leningrad, 1975.

\bibitem{Falkovich2002}
G.~Falkovich, A.~Fouxon, and M.~G. Stepanov.
\newblock Acceleration of rain initiation by cloud turbulence.
\newblock {\em Nature}, 419:151, 2002.

\bibitem{garcia1987monte}
Alejandro~L Garcia, Christian Van Den~Broeck, Marc Aertsens, and Roger
  Serneels.
\newblock A {M}onte {C}arlo simulation of coagulation.
\newblock {\em Physica A: Statistical Mechanics and its Applications},
  143(3):535--546, 1987.

\bibitem{liffman1992direct}
Kurt Liffman.
\newblock A direct simulation {M}onte-{C}arlo method for cluster coagulation.
\newblock {\em Journal of Computational Physics}, 100(1):116--127, 1992.

\bibitem{kruis2000direct}
F~Einar Kruis, Arkadi Maisels, and Heinz Fissan.
\newblock Direct simulation {M}onte {C}arlo method for particle coagulation and
  aggregation.
\newblock {\em AIChE Journal}, 46(9):1735--1742, 2000.

\bibitem{eibeck2000efficient}
Andreas Eibeck and Wolfgang Wagner.
\newblock An efficient stochastic algorithm for studying coagulation dynamics
  and gelation phenomena.
\newblock {\em SIAM Journal on Scientific Computing}, 22(3):802--821, 2000.

\bibitem{lee2000validity}
M.~H. Lee.
\newblock On the validity of the coagulation equation and the nature of runaway
  growth.
\newblock {\em Icarus}, 143(1):74--86, 2000.

\bibitem{MonteCarlo2003}
E.~Debry, B.~Sportisse, and B.~Jourdain.
\newblock A stochastic approach for the numerical simulation of the general
  dynamics equation for aerosols.
\newblock {\em Journal of Computational Physics}, 184(2):649--669, 2003.

\bibitem{matveev_tensor_2016}
S.~A. Matveev, D.~A. Zheltkov, E.~E. Tyrtyshnikov, and A.~P. Smirnov.
\newblock Tensor {T}rain versus {M}onte {C}arlo for the multicomponent
  {S}moluchowski coagulation equation.
\newblock {\em Journal of Computational Physics}, 316:164--179, 2016.

\bibitem{stadnichuk_smoluchowski_2015}
V.~Stadnichuk, A.~Bodrova, and N.~V. Brilliantov.
\newblock Smoluchowski aggregation–fragmentation equations: {Fast} numerical
  method to find steady-state solutions.
\newblock {\em International Journal of Modern Physics B}, 29(29):1550208,
  2015.

\bibitem{matveev_fast_2015}
S.~A. Matveev, A.~P. Smirnov, and E.~E. Tyrtyshnikov.
\newblock A fast numerical method for the {Cauchy} problem for the
  {Smoluchowski} equation.
\newblock {\em Journal of Computational Physics}, 282:23--32, 2015.

\bibitem{osinsky2020low}
AI~Osinsky.
\newblock Low-rank method for fast solution of generalized {S}moluchowski
  equations.
\newblock {\em Journal of Computational Physics}, 422:109764, 2020.

\bibitem{brilliantov2018increasing}
Nikolai~V Brilliantov, Arno Formella, and Thorsten P{\"o}schel.
\newblock Increasing temperature of cooling granular gases.
\newblock {\em Nature communications}, 9(1):1--9, 2018.

\bibitem{ZAGAYNOV2001983}
V.A. Zagaynov, K.~Denisenko, A.~Moskaev, and A.A. Lushnikov.
\newblock Periodical regimes in source-inhanced coagulating systems with sinks.
\newblock {\em Journal of Aerosol Science}, 32:983 -- 984, 2001.
\newblock Abstracts of the European Aerosol Conference 2001.

\bibitem{zagaynov2007periodic}
VA~Zagaynov, AA~Lushnikov, MS~Bakhtyreva, AO~Lutsenko, and TV~Khodzher.
\newblock Periodic regimes in the source-enhanced condensing aerodisperse
  system.
\newblock In {\em Doklady Earth Sciences}, volume 414, page 570. Springer
  Nature BV, 2007.

\bibitem{ball_collective_2012}
R.~C. Ball, C.~Connaughton, P.~P. Jones, R.~Rajesh, and O.~Zaboronski.
\newblock Collective {Oscillations} in {Irreversible} {Coagulation} {Driven} by
  {Monomer} {Inputs} and {Large}-{Cluster} {Outputs}.
\newblock {\em Physical Review Letters}, 109(16), October 2012.

\bibitem{krapivsky1994diffusion}
PL~Krapivsky.
\newblock Diffusion-limited-aggregation processes with three-particle
  elementary reactions.
\newblock {\em Physical Review E}, 49(4):3233, 1994.

\bibitem{matveev2019numerical}
SA~Matveev, DA~Stefonishin, AP~Smirnov, AA~Sorokin, and EE~Tyrtyshnikov.
\newblock Numerical studies of solutions for kinetic equations with
  many-particle collisions.
\newblock In {\em Journal of Physics: Conference Series}, volume 1163, page
  012008. IOP Publishing, 2019.

\bibitem{Bodrova2019}
A.~Bodrova, V.~Stadnichuk, P.~L. Krapivsky, J.~Schmidt, and N.~V. Brilliantov.
\newblock Kinetic regimes in aggregating systems with spontaneous and
  collisional fragmentation.
\newblock {\em J. Phys. A: Math. Gen.}, 52:205001, 2019.

\bibitem{pego2020temporal}
Robert~L Pego and Juan~JL Vel{\'a}zquez.
\newblock Temporal oscillations in {B}ecker--{D}{\"o}ring equations with
  atomization.
\newblock {\em Nonlinearity}, 33(4):1812, 2020.

\bibitem{niethammer2021oscillations}
Barbara Niethammer, Robert~L Pego, Andr{\'e} Schlichting, and Juan~JL
  Vel{\'a}zquez.
\newblock Oscillations in a becker-d$\backslash$" oring model with injection
  and depletion.
\newblock {\em arXiv preprint arXiv:2102.06751}, 2021.

\bibitem{KOB2017}
P.~L. Krapivsky, W.~Otieno, and N.~V. Brilliantov.
\newblock Phase transitions in systems with aggregation and shattering.
\newblock {\em Phys. Rev. E}, 96:042138, 2017.

\bibitem{matveev2017oscillations}
S.~A. Matveev, P.~L. Krapivsky, A.~P. Smirnov, E.~E. Tyrtyshnikov, and N.~V.
  Brilliantov.
\newblock Oscillations in aggregation-shattering processes.
\newblock {\em Physical Review Letters}, 119(26):260601, 2017.

\bibitem{timokhin2019newton}
IV~Timokhin, SA~Matveev, Nana Siddharth, Eugene~E Tyrtyshnikov, AP~Smirnov, and
  Nikolai~V Brilliantov.
\newblock Newton method for stationary and quasi-stationary problems for
  {S}moluchowski-type equations.
\newblock {\em Journal of Computational Physics}, 382:124--137, 2019.

\bibitem{French}
R.~S. French, S.~K. Hicks, M.~R. Showalter, A.~K. Antonsen, and D.~R. Packard.
\newblock Analysis of clumps in saturn’s f ring from voyager and cassini.
\newblock {\em Icarus}, 241:200, 2014.

\bibitem{budzinskiy2020hopf}
Stanislav~S Budzinskiy, Sergey~A Matveev, and Pavel~L Krapivsky.
\newblock Hopf bifurcation in addition-shattering kinetics.
\newblock {\em arXiv preprint arXiv:2012.09003}, 2020.

\bibitem{gillespie1975exact}
Daniel~T Gillespie.
\newblock An exact method for numerically simulating the stochastic coalescence
  process in a cloud.
\newblock {\em Journal of the Atmospheric Sciences}, 32(10):1977--1989, 1975.

\bibitem{matsoukas1997monte}
Themis Matsoukas and Y~Tang.
\newblock Monte {C}arlo simulation of agglomeration and grinding.
\newblock {\em Particulate Science and Technology}, 15(2):156--156, 1997.

\bibitem{xu2014fast}
Zuwei Xu, Haibo Zhao, and Chuguang Zheng.
\newblock Fast {M}onte {C}arlo simulation for particle coagulation in
  population balance.
\newblock {\em Journal of aerosol science}, 74:11--25, 2014.

\bibitem{wei2013gpu}
J~Wei and F~Einar Kruis.
\newblock G{P}{U}-accelerated {M}onte {C}arlo simulation of particle
  coagulation based on the inverse method.
\newblock {\em Journal of Computational Physics}, 249:67--79, 2013.

\bibitem{xu2015accelerating}
Zuwei Xu, Haibo Zhao, and Chuguang Zheng.
\newblock Accelerating population balance-{M}onte {C}arlo simulation for
  coagulation dynamics from the {M}arkov jump model, stochastic algorithm and
  {G}{P}{U} parallel computing.
\newblock {\em Journal of Computational Physics}, 281:844--863, 2015.

\bibitem{BrillPRE2018}
N.~V. Brilliantov, W.~Otieno, S.~A. Matveev, A.~P. Smirnov, E.~E. Tyrtyshnikov,
  and Krapivsky~P. L.
\newblock Steady oscillations in aggregation-fragmentation processes.
\newblock {\em Phys. Rev. E}, 98:012109, 2018.

\bibitem{zacharov2019zhores}
Igor Zacharov, Rinat Arslanov, Maksim Gunin, Daniil Stefonishin, Andrey Bykov,
  Sergey Pavlov, Oleg Panarin, Anton Maliutin, Sergey Rykovanov, and Maxim
  Fedorov.
\newblock “{Z}hores”—{P}etaflops supercomputer for data-driven modeling,
  machine learning and artificial intelligence installed in {S}kolkovo
  {I}nstitute of {S}cience and {T}echnology.
\newblock {\em Open Engineering}, 9(1):512--520, 2019.

\end{thebibliography}

\end{document}